\documentclass[12pt,reqno]{amsart}
\usepackage{amsmath,amssymb}
\usepackage{graphicx,xcolor}
\usepackage{latexsym}
\usepackage{mathrsfs}
\usepackage[abbrev]{amsrefs}
 
\newtheorem{thm}{Theorem}[section]
\newtheorem{cor}[thm]{Corollary}
\newtheorem{prop}[thm]{Proposition}
\newtheorem{lem}[thm]{Lemma}

\newtheorem{definition}[thm]{Definition}

\newcommand{\R}{\mathbb R}

\newcommand{\BLACK}{\color{black}}
\newcommand{\RED}{\color{red}}

 \newenvironment{pf}
    {{\noindent \bf Proof. }}{\hfill $\Box$}

\numberwithin{equation}{section}
\numberwithin{thm}{section}

\begin{document}

\begin{center}
\large \bf 
Sobolev spaces on arbitrary domains and semigroups 
\\
generated by fractional Laplacian
\end{center}

\footnote[0]
{
{\it Mathematics Subject Classification} (2010): 
Primary 46E35; 
Secondary 46F12.

{\it Keywords}: 
Sobolev spaces, arbitrary domains, semigroup, fractional Dirichlet Laplacian, embeddings, Gagliardo-Nirenberg inequality

{\it Addresses}: 
Reinhard Farwig, Fachbereich Mathematik, Technische, Universit\"at Darmstadt, 64283 Darmstadt, Germany, 
{\tt farwig@mathematik.tu-darmstadt.de}

Tsukasa Iwabuchi, Mathematical Institute, Tohoku University, 
Sendai, 980-8578, Japan, 
{\tt t-iwabuchi@tohoku.ac.jp}
}
\vskip5mm

\begin{center}
Reinhard FARWIG \quad and \quad 
Tsukasa IWABUCHI  

\vskip2mm

\end{center}

\vskip5mm

\begin{center}
\begin{minipage}{135mm}
\footnotesize
{\sc Abstract. } 
We describe a procedure to introduce Sobolev spaces and the semigroup generated by the fractional 
Dirichlet Laplacian on an arbitrary domain of $\R^d$. In particular, the well-definedness of the spaces 
of both non-homogeneous and homogeneous type together with their duality properties, embeddings, and Gagliardo-Nirenberg inequalities will be discussed.
We also show the continuity and the smoothing property of the semigroup. 

\end{minipage}
\end{center}

\section{Introduction} 

We study the theory of Sobolev spaces on general domains. 
The spaces that are nowadays called the Sobolev spaces 
$W^{m,p}$, with $m = 0,1,2,\ldots$, were introduced by Sobolev~\cites{Sobo-1935,Sobo-1936,Sobo-1938} 
in 1935--1938. 
Related to them are the Bessel potential spaces $H^s_p$, $s\in\R$, introduced by 
Aronszajn-Smith~\cite{AS-1961} and Calder\'on~\cite{Ca-1961} in 1961, 
which are spaces with fractional derivatives on the entire Euclidean space 
such that $H^m_p (\R^d) = W^{m,p} (\mathbb R^d)$, $1 < p < \infty$. 
Moreover, $H^s_p$ is a complex interpolation space between $ W^{m,p}$ and $W^{m+1,p} $ for $m<s<m+1$.
The aim of this paper is to define such kind of spaces of 
non-homogeneous and also homogeneous type on general domains 
and to investigate properties of the semigroup generated by the fractional Laplacian, 
motivated by the applications to partial differential equations.

There are several ways for defining Sobolev spaces. 
The original definition is based on the integrability of functions 
and their weak derivatives, an idea which applies to any domain; 
but it is not obvious how to introduce the fractional Laplacian and semigroups on 
general domains. 
On the entire space $\mathbb R^d$, the Fourier transform is important to 
characterize not only Sobolev spaces but also 
Besov spaces and other spaces. However, such a useful tool is not available on domains. 
If the boundary of a domain is smooth, we can consider the zero extension of functions 
on the domain to the entire space and, {\em vice versa},  
the restriction of functions on $\mathbb R^d$ to the domain. 
We also mention that norms on spaces with non-integer exponent of regularity can be defined by double integrals of difference quotients. 
To the best of our knowledge, we need some restriction on domains or indices 
and we refer to \cites{DiPaVa-2012,Triebel_1983,Triebel_1995} for these characterizations 
and the history of function spaces. See \cite{Ryc-1999} for corresponding results on Lipschitz domains.

In this paper, we define Sobolev spaces of non-homogeneous and homogeneous type, $H^{s,p}$ and $\dot H^{s,p}$, respectively, on an arbitrary domain without 
any restriction on indices $s \in \mathbb R, 1 \leq p \leq \infty$ related to the Dirichlet Laplacian. In addition, we study their properties 
such as duality, Sobolev embeddings, Gagliardo-Nirenberg inequalities, {\em etc.}, and also analyze
the semigroups generated by the Dirichlet Laplacian of fractional order. 
It should be noted that our theory can be based also on \BLACK more general operators satisfying Gaussian upper bounds of the kernel, $K$, of the semigroup, {\em i.e.}
$0\leq K(x,y,t) \leq Ct^{-n/2} \exp(-{c|x-y|^2}/{t})$. 

Let $\Omega$ be an arbitrary domain in $\mathbb R^d $ with $d \geq 1$. 
We denote by the Dirichlet 
Laplacian, $A$, the operator defined by
\begin{equation}\notag 
\begin{cases}\mathcal D (A) 
= \big\{ f \in H^1_0(\Omega)  \, \big| \, 
     \Delta f \in L^2 (\Omega) \big\},
\\
\displaystyle 
Af = -\Delta f 
= - \sum _{j=1}^d \frac{\partial^2}{\partial x_ j^2 } f, 
\quad f \in \mathcal D(A). 
\end{cases}
\end{equation}
Actually, $A$ is a non-negative self-adjoint operator given by its form definition such that the domain of $\sqrt{A} = A^{1/2}$ equals $H^1_0(\Omega)$. 
Let $\{ E(\lambda) \}_{\lambda}$ denote the resolution of the identity given by $A$, introduced by the spectral theorem and supported on $[0,\infty)$. 
The operator $A$, defined initially on $H^1_0(\Omega)$, is also regarded as an operator 
on spaces of distributions $\mathcal X'(A)$, $\mathcal Z'(A)$ defined below.

Let $\phi_0 \in C^\infty_0(\mathbb R)$ be a non-negative function on $\mathbb R$ such that  
\begin{equation}
\label{917-1}
{\rm supp \, } \phi _0
\subset \{ \, \lambda \in \mathbb R \, | \, 2^{-1} \leq \lambda \leq 2 \, \}, 
\quad \sum _{ j \in \mathbb Z} \phi_0 ( 2^{-j}\lambda) 
 = 1 
 \quad \text{for } \lambda > 0,  
\end{equation}
and let $ \{ \phi_j \}_{j \in \mathbb Z}$ be defined by  
\begin{equation} \label{917-2}
\phi_j (\lambda) = \phi_0 (2^{-j} \lambda) \quad \text{for } \lambda \in \mathbb R . 
\end{equation}

\vskip1mm

\noindent 
\begin{definition}\label{D1.1}
(i) 
The test function space $\mathcal X (A)$ {\em of non-homogeneous type} is 
defined by 
\begin{equation}\notag 
\mathcal X (A)
:= \big\{ f \in  L^1 (\Omega) \cap \mathcal D (A) 
 \, \big| \, 
    A^{m} f \in L^1(\Omega ) \cap \mathcal D (A) \text{ for all } m \in \mathbb N 
   \big\} 
\end{equation} 
equipped with the family of semi-norms $\{ p_{m} (\cdot) \}_{ m = 1 } ^\infty$ 
given by 
\begin{equation}\notag 
p_{m}(f) := 
\| f \|_{ L^1(\Omega)} 
+ \sup _{j \in \mathbb N} 2^{mj} 
  \| \phi_j (\sqrt{A}) f \|_{ L^1(\Omega)} . 
\end{equation}
Then let $\mathcal X'(A)$ denote the topological dual of $\mathcal X (A)$.

(ii) The test function space $\mathcal Z(A)$ {\em of homogeneous type} is defined by 
\begin{equation}\notag 
\mathcal Z (A) 
:= \Big\{ f \in \mathcal X (A) 
 \, \Big| \, 
  \sup_{j \leq 0} 2^{ m |j|} 
    \big\| \phi_j \big(\sqrt{ A } \big ) f \big \|_{L^1(\Omega)} < \infty 
  \text{ for all } m \in \mathbb N
   \Big\} 
\end{equation}
equipped 
with the family of semi-norms $\{ q_{m} (\cdot) \}_{ m = 1}^\infty$ given by 
\begin{equation}\label{EQ:qV} 
q_{m}(f) := 
\| f \|_{L^1 (\Omega) }
+ \sup_{j \in \mathbb Z} 2^{m|j|} \| \phi_j (\sqrt{A}) f \|_{L^1(\Omega)}. 
\end{equation}
Moreover, $\mathcal Z' (A)$ denotes the topological dual of $\mathcal Z (A)$. 
\end{definition}

It should be noted that the operator $\phi_j (\sqrt{A})$ on $L^1 (\Omega) $ 
is initially defined by the 
spectral decomposition 
\[
\phi_j(\sqrt{A})  = \int _0^\infty  \phi_j(\sqrt{\lambda})\,{\rm d} E(\lambda) 
\quad \text{ on } L^2(\Omega) 
\]
and then extended to $L^1(\Omega)$ by the spectral multiplier theorem 
(see \cites{IMT-2018,TOS-2002} and also Lemma~\ref{lem:1} below). 
We see that $\mathcal X(A)$, $\mathcal Z(A)$ are Fr\'echet spaces independent of the choice of 
$\{ \phi_j \}_{j \in \mathbb Z}$ and their duals can be regarded as spaces of distributions 
associated with the Dirichet Laplacian. 
As a comparison with the whole space, if $\Omega = \mathbb R^d$, then 
$\mathcal S(\mathbb R^d) \hookrightarrow \mathcal X(A) 
$ and 
$
\mathscr Z(\mathbb R^d) \hookrightarrow \mathcal Z(A) 
$. 
Here 
$\mathcal S(\mathbb R^d)$ is the Schwartz class and 
\[
\mathscr Z(\mathbb R^d) = 
\Big\{ f \in \mathcal S(\mathbb R^d) \, | \, 
  \int_{\mathbb R^d} x^\alpha f(x) \,{\rm d}x  = 0 
  \text{ for all } \alpha \in (\mathbb N_0)^d 
\Big\},
\]
{\em cf.} \cite[\S 5.1.2]{Triebel_1983}. 
Moreover, for $f \in \mathcal S(\mathbb R^d)$, we see that
$f \in \mathscr Z(\mathbb R^d)$ if and only if $f \in \mathcal Z(A)$. 
From that point of view \BLACK the condition on vanishing moments in $\mathscr Z(\mathbb R^d)$ is replaced in the space $\mathcal Z(A)$ by the condition that $\sup_{j\leq 0} 2^{m|j|} 
    \big\|\phi_j \big(\sqrt{A} \big ) f \big \|_{L^1(\Omega)} < \infty$.
We here recall the relation between $\mathcal X(A), \mathcal Z(A)$, 
the Lebesgue spaces, and their duals. 

\begin{lem}\cite[Lemma~4.6]{IMT-2019} \label{lem:0301-3}
The following inclusions hold for all $1 \leq p \leq \infty$. 
\begin{align*}
\mathcal X(A), \mathcal Z(A) & \subset L^p (\Omega) \subset \mathcal X'(A), \mathcal Z'(A).
\end{align*}
\end{lem}
\BLACK

Next we define the Sobolev spaces we are looking for as subspaces of $\mathcal X'(A)$ and $\mathcal Z'(A)$.

\vskip1mm

\noindent 
\begin{definition}\label{D1.2} 
Let $s \in \mathbb R$ and $1 \leq p \leq \infty$. 
\begin{enumerate}
\item[(i)] 
The Sobolev spaces of non-homogeneous type are defined by 
\[
H^s_p (A) := 
\left\{ f \in \mathcal X'(A) \, | \, 
   \| f \|_{H^s_p (A)} := \left\| (1+A)^{\frac{s}{2}} f  \right\|_{L^p (\Omega) } < \infty
\right\}. 
\]

\item[(ii)] 
The Sobolev spaces of homogeneous type are defined by 
\[
\dot H^s_p (A) := 
\left\{ f \in \mathcal Z'(A) \, | \, 
 \| f \| _{\dot H^s_p} := \| A^{\frac{s}{2}} f \|_{L^p (\Omega)} < \infty
\right\}  .
\]
\end{enumerate}
\end{definition}

We state three theorems on the well-definedness of these Sobolev spaces, 
their properties and the semigroup generated by fractional powers of $A$. 
The first theorem concerns well-definedness of the spaces, duality and embeddings. 

\begin{thm}\label{thm:1}
Let $s \in \mathbb R$, $1 \leq p \leq \infty$, and let $p'\in[1,\infty]$ be the conjugate exponent to $p$.

\begin{enumerate}
\item[(i)] 
$H^s_p (A)$ and $\dot H^s_p(A)$ are Banach spaces and they satisfy the embeddings
\[
\mathcal X(A) \hookrightarrow H^s_p (A) \hookrightarrow \mathcal X'(A)  , 
\quad 
\mathcal Z(A) \hookrightarrow \dot H^s_p (A) \hookrightarrow \mathcal Z'(A) .
\]

\item[(ii)] 
{\rm (}Dual spaces{\rm )}
If $1 \leq p < \infty$, then 
$$ (H^s_p(A))' \cong H^{-s}_{p'}(A) \quad\textrm{ and}\quad  (\dot H^s_{p}(A))'\cong \dot H^{-s}_{p'}(A). $$ 

\item[(iii)] 
{\rm (}Lift property{\rm )}
Let $s _0 \in \mathbb R$. Then the maps
\begin{align*} 
(1+A)^{s_0 /2} : H^{s}_p (A) \to H^{s-s_0}_p (A),\; & \;f\mapsto (1+A)^{s_0 /2}f,\\   
A^{s_0 /2} : \dot H^{s}_p (A) \to \dot H^{s-s_0}_p (A),\; & \;f\mapsto A^{s_0 /2}f 
\end{align*} 
are isomorphisms.  

\item[(iv)] 
{\rm (}Sobolev embedding{\rm )}
Let $1 < r \leq p < \infty$. Then 
\[
H^{s+d(\frac{1}{r}- \frac{1}{p})}_r (A) \hookrightarrow H^s_p(A), 
\quad 
\dot H^{s+d(\frac{1}{r}- \frac{1}{p})}_r (A) \hookrightarrow \dot H^s_p(A). 
\]

\item[(v)] 
Let $s < d/p$. Then 
\[
\dot H^s_p (A) 
\simeq 
\Big\{ f \in \mathcal X'(A) \, \Big| \, 
\| f \|_{\dot H^s_p (A)} < \infty , \,\,
 f = \sum _{j \in \mathbb Z} \phi_j (\sqrt{A}) f \text{ in } \mathcal X'(A)
\Big\}. 
\]
\end{enumerate}
\end{thm}

The second theorem studies the Gagliardo-Nirenberg interpolation inequality. 
The original was proved by Gagliardo~\cite{Gag-1959} 
and Nirenberg~\cite{Nir-1959} independently. 
We refer to Hajaiej {\em et al.} \cite{HMOW-2011} showing the necessity and 
the sufficiency of indices obtained in the whole space, 
and to Chikami \cite{Chi-2018} for the inequality in the Fourier-Herz spaces; see
recent papers by Brezis \& Mironescu  \cites{BrMi-2018,BrMi-2019}  investigating the validity of the inequality in a standard domain. 

The following result holds on domains without any assumption on the boundary. 

\begin{thm}\label{thm:2}
 {\rm(}Gagliardo-Nirenberg inequality{\rm)}
Let $s,s_0 > 0$, $ 1 \leq p,r,r_0 \leq \infty$, and
$\theta \in (0,1)$ satisfy  
\begin{equation}\label{0912-1}
s- \frac{d}{p} 
= \theta \Big( - \frac{d}{r} \Big) + (1-\theta) \Big( s_0 - \frac{d}{r_0} \Big) ,
\end{equation}
\begin{equation}\notag 
-\frac{d}{r} \not = s_0 - \frac{d}{r_0}, 
\qquad 
\begin{cases}
s \leq (1-\theta) s_0 
& \text{if } \max\{r, r_0\} \leq p,
\\
s < (1-\theta) s_0
& \text{if } \min\{r, r_0 \} \leq p  < \max\{r, r_0\}, \BLACK
\end{cases}
\end{equation}
Then 
\[
\| f \|_{\dot H^s_p (A)} 
\leq C \| f \|_{L^{r}} ^\theta \| f \|_{\dot H^{s_0}_{r_0}} ^{1-\theta}. 
\]
\end{thm}

\noindent 
{\bf Remark. } 
The excluded case when 
$s= (1-\theta) s_0 $ and $ \max\{ r, r_0 \} > p $ 
is studied in the whole space case (see \cite{HMOW-2011}). 
We conjecture an affirmative answer if 
$1 < p ,r,r_0 < \infty$. 

\vskip3mm

The third 
theorem is concerned with the semigroup generated by the 
Dirichlet Laplacian of fractional order, $A^{\alpha/2}.$ 
We define the fractional Laplacian and the semigroup $\{e^{-tA^{\alpha/2}}\}_{t\geq 0}$ via the spectral multipliers $\phi_j(\sqrt{A})$ and by duality, see Definitions~\ref{def:1}, \ref{def:2} and subsequent propositions.

The following results for the semigroup $\{e^{-tA^{\alpha/2}}\}$ hold on Sobolev spaces of homogeneous and non-homogeneous type; they were investigated for Besov spaces of homogeneous type in \cite{Iwa-2018}. 

\begin{thm}\label{thm:3} 
{\rm (}Semigroups in $\dot H^s_p(A), H^s_p(A)${\rm )} 
Let $\alpha > 0$, $s \in \mathbb R$, $1 \le p \leq \infty $ and $t >0$. 
\begin{enumerate}
\item[(i)] 
$e^{-t A^{\frac{\alpha}{2}}}$ is a bounded linear operator in $\dot H^s_p (A)$, 
{\em i.e.}, 
there exists a constant $C>0$ independent of $t$ such that 
for any $ f \in \dot H^s_p(A)$ 
\[
e^{-tA^{\frac{\alpha}{2}}} f \in \dot H^s_p (A) 
\quad \text{and} \quad 
\big\| e^{-tA^{\frac{\alpha}{2}}} f  \big\|_{\dot H^s_p (A)} \leq C \|f\|_{\dot H^s_p (A)}. 
\]
If the domain is bounded, then $e^{-t A^{\frac{\alpha}{2}}}$ decays exponentially, 
{\em i.e.,} there exists $\kappa > 0$ such that 
\[
\big\| e^{-tA^{\frac{\alpha}{2}}} f  \big\|_{\dot H^s_p (A)}  
\leq C e^{-\kappa t} \|f\|_{\dot H^s_p (A)} 
\quad \text { for all } t > 0 . 
\]
\BLACK 
The same result holds for the non-homogeneous spaces $H^s_p(A)$.

\item[(ii)] 
Let 
\begin{equation}\label{0807-1}
s_2 \geq s_1 , \quad p_2 \geq p_1 , \quad 
d \Big(\frac{1}{p_1} - \frac{1}{p_2} \Big) + s_2 - s_1 > 0. 
\end{equation}
Then there exists a constant $C > 0$ independent of $t$ such that 
\begin{equation}\notag 
\big\| e^{-tA^{ \frac{\alpha }{2} }} f \big\|_{\dot H^{s_2}_{p_2} (A)} 
\leq C 
t^{- \frac{s_2 - s_1}{\alpha}-\frac{d}{\alpha} (\frac{1}{p_1} - \frac{1}{p_2}) 
         }      \| f \|_{ \dot H^{s_1}_{p_1} (A) } 
\end{equation}
for any $f \in \dot H^{s_1}_{p_1} (A)$. 
In the non-homogeneous spaces $H^s_p(A)$, the difference in regularity, 
$s_2 -s_1$, is effective only for small $t$ and, for any $f \in  H^{s_1}_{p_1} (A)$, \BLACK
\begin{equation}\notag 
\big\| e^{-tA^{ \frac{\alpha }{2} }} f \big\|_{H^{s_2}_{p_2} (A)} 
\leq C 
\Big( 1 + t^{- \frac{s_2 - s_1}{\alpha}} \Big) 
t^{-\frac{d}{\alpha} (\frac{1}{p_1} - \frac{1}{p_2}) 
         }      \| f \|_{ H^{s_1}_{p_1} (A) }  .
\end{equation}
If the domain is bounded, then in both estimates an additional exponentially decaying term as in (i) can be inserted. 
\item[(iii)] 
Assume that $1 \leq p < \infty$ and $f \in \dot H^s_p (A), H^s_p (A)$. 
Then 
\[
\lim_{ t \to 0} 
 \big\| e^{-t A^{\frac{\alpha}{2}}} f - f \big\|_{\dot H^s_p(A)} = 0 , 
 \quad \lim_{ t \to 0} 
 \big\| e^{-t A^{\frac{\alpha}{2}}} f - f \big\|_{ H^s_p(A)} = 0 , 
\]
respectively. 
If $p = \infty$, then $e^{-t A^{\frac{\alpha}{2}}} f$ converges to $f$ 
in the weak-$*$ sense, {\em i.e.} 
\[
\lim_{ t\to 0} 
\int_{\Omega} 
\Big( A^{s/2}\big( e^{-t A^{\frac{\alpha}{2}}} f - f \big) \Big)  \overline{A^{-s/2} g} \,{\rm d}x = 0 
\quad \text{for any } g \in \dot H^{-s}_{1} (A), 
\]
\[
\lim_{ t\to 0} 
\int_{\Omega} 
\Big( (1+A)^{s/2}\big( e^{-t A^{\frac{\alpha}{2}}} f - f \big) \Big)  \overline{(1+A)^{-s/2} g} \,{\rm d}x = 0
\quad \text{for any } g \in  H^{-s}_{1} (A). 
\]

\item[(iv)]  
For $1\leq p\leq\infty$ the mapping $t \mapsto e^{-t A^{\frac{\alpha}{2}}} $ is analytic 
from $(0,\infty) \to \mathcal L(\dot H^s_p (A))$. If $1\leq p<\infty$, then $-A^{\frac{\alpha}{2}}$ generates an analytic semigroup. 
The same results hold for the non-homogeneous spaces $H^s_p(A)$. 
\end{enumerate}
\end{thm}

\vskip3mm 

\noindent 
{\bf Remark. } 
As a corollary of (i), (ii) in Theorem~\ref{thm:3}, 
the following inclusions are dense:  
\[
\dot H^s_p (A) \cap \dot H^{s+\delta}_p (A) 
\subset \dot H^s_p (A),
\quad 
H^{s+\delta}_p (A) 
\subset H^s_p (A), 
\quad \delta > 0 , 1 \leq p < \infty .
\]
In (iv) of Theorem~\ref{thm:3}, 
each $A^{\frac{\alpha}{2}}$  is considered as an operator 
$\dot H^s_p(A) \cap \dot H^{s+\alpha}_p(A) \to \dot H^s_p(A) $ and 
$
H^{s+\alpha}_p(A) \to H^s_p(A) $, 
respectively, with dense domain when $1 \leq p < \infty$. 
\BLACK

\vskip3mm

Let us give some comments for the proof of the theorems. 
The proof of Theorem~\ref{thm:1} is similar to the paper \cite{IMT-2019} except for the 
Sobolev embedding, which is shown by a proof of decay estimates of the heat kernel 
implying the embedding (see {\em e.g.} Davies~\cite{Davies_1990}). 
Theorem~\ref{thm:2} is obtained by proving a similar result in Besov spaces analogously to the whole space case 
(see {\em e.g.} \cites{HMOW-2011,MaOz-2003}) together with the spectral 
multiplier theorem (Lemma~\ref{lem:1}) instead of the Fourier multiplier theorem.
We can prove a part of Theorem~\ref{thm:3} 
(i), (ii), (iii) for 
the high spectral component of functions, analogously to the paper \cite{Iwa-2018}. 
As for the low spectral part, we can not apply the spectral multiplier theorem 
because of little regularity of the functions 
$e^{-t \lambda ^{\alpha /2}} $, $0 < \alpha <1$, 
which causes a difficulty in the non-homogeneous case. 
Instead of the boundedness of the spectral multipliers, 
we apply the subordination principle (see the monograph by Yosida~\cite[Chapter IX.11]{Yosida_1980}) which can be understood by using the formula:
\begin{equation}\notag 
e^{-t \lambda ^\alpha} = \int_0^\infty F_{t,\alpha} (s) e^{-s\lambda } \,{\rm d}s , 
\quad \lambda > 0 , 0 < \alpha < 1,
\end{equation}
where $F_{t,\alpha}$ is called the subordinator (see \eqref{0829-1} below). 
This allows to estimate the non-smooth function $e^{-t \lambda ^\alpha}$ 
by the smooth one, $e^{-s\lambda }$. 
The analyticity (iv) of the semigroup follows from 
the fact that resolvent estimates imply sectoriality and analyticity 
(see \cite{Lunardi_1995}), 
where resolvent estimates will be obtained by spectral multiplier estimates. 
We also mention that analyticity of the semigroup in the fractional case can be based on the non-fractional case 
(see \cite[page 263, Theorem 1]{Yosida_1980}). 

\vskip3mm 

\noindent 
{\bf Remark. } 
Since the low spectral component is handled by subordination, 
we obtain the result for the non-homogeneous Besov spaces 
corresponding to Theorem~\ref{thm:3}, which was not treated in the paper~\cite{Iwa-2018}. 

\vskip3mm 

This paper is organized as follows. 
In Section 2, we recall the spectral multiplier theorem and define the
fractional Laplacian and the semigroups generated by them using the subordination principle. 
Sections 3, 4 and 5 are devoted to proving Theorems \ref{thm:1} -- \ref{thm:3}.

\section{Preliminaries}

\begin{lem}\label{lem:1} 
{\rm  (Spectral multiplier estimates \cites{IMT-2018,TOS-2002})} 
Let $\varphi \in\mathcal S(\mathbb R)$, $1 \leq p \leq \infty$, 
$\beta > (d+1)/2 \, (\beta \in \mathbb N)$ and $\gamma  > (3d+2)/4 $.  
Then 
\begin{equation}\label{0806-1}
 \| \varphi (t A) \|_{L^p \to L^p} 
\leq C \| (1+|\cdot|^2)^{\frac{\gamma}{2}} \varphi \|_{H^{\beta}} 
\text{ for all } t > 0 .
\end{equation}
Moreover, if $s \in \mathbb R$, $1 \leq r \leq p \leq \infty$ and 
$\widetilde\gamma > (3d+2)/4 + d(1/r - 1/p)/2$,  
then 
\begin{equation}\label{0802-1}
\| A^{s/2} \phi_j (\sqrt{A}) \|_{L^r \to L^p} 
\leq C 2^{sj + d(\frac{1}{r} - \frac{1}{p}) j} 
\| (1+|\cdot|^2)^{\frac{\widetilde\gamma}{2}} |\cdot |^s \phi_0 \|_{H^{\beta}} 
\quad \text{for all } j \in \mathbb Z .
\end{equation}
\end{lem}
%

By Lemma~\ref{lem:1}, 
we can consider the operator $\varphi (A)$ on $L^p$ 
provided that $\varphi $ belongs to $H^\beta$ with compact support.  
For the proof of the boundedness of the spectral multipliers we refer to  
\cite[Theorem~1.1, Section 8]{IMT-2018}, \cite[Theorem~1.1]{IMT-pre} 
and \cite[Proposition~6.1]{TOS-2002}. 

We recall the resolution of the identity.  
Let $\psi \in C_0 ^\infty (\mathbb R)$ be such that 
\begin{equation}\label{0301-1}
{\rm supp } \, \psi \subset \Big[ -\frac{1}{2} , 2\Big], 
\qquad \psi (\lambda) + \sum _{ j =1}^\infty \phi_j (\lambda) = 1
\quad \text{for all } \lambda > 0 .
\end{equation}

\begin{lem}\cite[Lemma~4.5]{IMT-2019} \label{lem:5} 
{\rm (i)} Let $f \in \mathcal X(A)$ or $\mathcal X'(A)$. Then 
\[
f = \psi (\sqrt{A}) f +  \sum_{j=1}^\infty \phi_j (\sqrt{A}) f
\quad \text{ in }\; \mathcal X(A) \; \text{ or }\;   \mathcal X'(A),   \text{ respectively}.
\]
where in the latter case 
$\psi (\sqrt{A}) f$ and $\phi_j (\sqrt{A})f$ are considered as functions in $L^\infty(\Omega)$.  
\\
{\rm (ii)} 
Let $f \in \mathcal Z(A)$ or $\mathcal Z'(A)$. Then 
\[
f = \sum _{j \in \mathbb Z} \phi_j (\sqrt{A}) f 
\quad \text{ in }\; \mathcal Z(A) \; \text{ or }\;   \mathcal Z'(A),   \text{ respectively},\] 
where in the latter case 
$\psi (\sqrt{A}) f$ and $\phi_j (\sqrt{A})f$ are considered 
as functions in $L^\infty(\Omega)$. 
\end{lem}

\vskip3mm 

We define the fractional Laplacian and the semigroup in the following way.

\vskip3mm 

\begin{definition}\label{def:1}
(Fractional Laplacian)
\begin{enumerate}
\item[(i)] 
Let $s > 0$. 
The operator $A^{s/2}$ on $\mathcal X(A)$ is defined by 
\[
A^{s /2} f := \sum _{j \in \mathbb Z} A^{s/2} \phi_j (\sqrt{A}) f 
\quad \text{in } \mathcal X(A) , 
\quad f \in \mathcal X(A). 
\]
The operator $A^{s/2}$ on $\mathcal X'(A)$ is defined by 
\[
_{\mathcal X'}\langle  A^{s/2} f , g \rangle _{\mathcal X} 
:= \,  _{\mathcal X'}\langle  f , A^{s/2} g \rangle _{\mathcal X} , 
\quad f \in \mathcal X'(A),  g \in \mathcal X(A) .
\]

\item[(ii)]  
Let $s \in \mathbb R$. 
The operator $(1+A)^{s/2}$ on $\mathcal X(A)$ is defined by 
\[
(1+A)^{s/2} f := \sum _{j \in \mathbb Z}(1+A)^{s/2} \phi_j (\sqrt{A}) f 
\quad \text{in } \mathcal X(A) , 
\quad f \in \mathcal X(A). 
\]
The operator $(1+A)^{s/2}$ on $\mathcal X'(A)$ is defined by 
\[
_{\mathcal X'}\langle  (1+A)^{s/2} f , g \rangle _{\mathcal X} 
:= \,  _{\mathcal X'}\langle  f , (1+A)^{s/2} g \rangle _{\mathcal X} , 
\quad f \in \mathcal X'(A),  g \in \mathcal X(A) .
\]
 
\item[(iii)] 
Let $s \in \mathbb R$. 
The operator $A^{s/2}$ on $\mathcal Z(A)$ is defined by 
\[
A^{s /2} f := \sum _{j \in \mathbb Z} A^{s/2} \phi_j (\sqrt{A}) f 
\quad \text{in } \mathcal Z(A) , 
\quad f \in \mathcal Z(A). 
\]
Moreover, the operator $A^{s /2}$ on $\mathcal Z'(A)$ is defined by 
\[
_{\mathcal Z'}\langle  A^{s/2} f , g \rangle _{\mathcal Z} 
:= \,_{\mathcal Z'}\langle  f , A^{s/2} g \rangle _{\mathcal Z} , 
\quad f \in \mathcal Z'(A), g \in \mathcal Z(A).
\]
\end{enumerate}
\end{definition}

\begin{prop}
Let $s >0$. The operators $A^{s/2}$ on $\mathcal X(A)$, $\mathcal X'(A)$, 
$\mathcal Z(A)$, and $\mathcal Z'(A)$ are well-defined and continuous. 
\end{prop}

\begin{pf}
We start with the non-homogeneous case. 
To show the well-definedness, it suffices to verify that 
$A^{s/2}$ on $\mathcal X(A)$ is well-defined and continuous, since 
$A^{s/2}$ on $\mathcal X'(A)$ is defined by duality.
Let $f \in \mathcal X(A)\subset L^2$. 
We use the resolution of the identity on $L^2$, take into account that zero is not an eigenvalue of $A$
(see also the generalized decomposition of $\mathcal X(A)$ in Lemma~\ref{lem:5}) and conclude that 
\[
f (x)
= \sum _{j=-\infty}^\infty \phi_j( \sqrt{A}) f  (x)
\quad \text{for almost every } x \in \Omega.
\]
We have from \eqref{0802-1} that 
\[
\sum _{ j \leq 0} p_m \Big(  A^{s /2}\phi_j( \sqrt{A}) f \Big)  
\leq C \Big(\sum_{j \leq 0} 2^{s j} \Big) \| f \|_{L^1 (\Omega)} ,
\]
and for $m_s \in \mathbb N$ with $m_s > m + s $ 
\[
\sum _{ j =1}^\infty p_m \Big( A^{s /2}\phi_j( \sqrt{A}) f \Big) 
\leq C \sum _{j=1}^\infty 2^{ (s + m) j } \| \phi_j (\sqrt{A}) f\|_{L^1 (\Omega)}
\leq C p_{m_s} (f) . 
\]
These estimates naturally lead to the justification of the definition
\[
A^{s/2}f 
= \sum _{j=-\infty}^\infty A^{s/2} \phi_j( \sqrt{A}) f  
\quad \text{ in } \mathcal X(A)  , 
\quad f \in \mathcal X(A),
\]
and the series in the right hand side converges in $\mathcal X(A)$. 
The continuity of $A^{s/2}$ on $\mathcal X(A)$ holds by applying the estimates above. 
Thus, the well-definedness of $A^{s/2}$ on $\mathcal X'(A)$ is also obtained, 
since the operator is defined by the dual operator. 

The homogeneous case follows analogously. 
\end{pf}

\vskip3mm

For a negative order, 
the well-definedness of $A^{-s /2}$, $s > 0$, on $\mathcal Z'(A)$ can be obtained 
by estimating the low spectral component 
\[
\sum _{ j \leq 0}  q_m \Big(  A^{-s /2}\phi_j( \sqrt{A}) f \Big)  
\leq C q_{m_s } (f) ,
\]
where $m_s \in\mathbb N$ satisfies $m_s > m+ s$.  
In fact, it follows by Lemma~\ref{lem:1} that 
\[
\begin{split}
q_m \Big(  A^{-s /2}\phi_j( \sqrt{A}) f \Big) 
\leq 
& C 2^{-sj} q_m \Big( \phi_j( \sqrt{A}) f \Big)
\leq C  q_{ m+s} \Big( \phi_j( \sqrt{A}) f \Big) 
\\
\leq
&  C 2^{(-m_s + m + s)|j|} q_{m_s} (f) .
\end{split}
\]
which proves the inequality above by taking the sum over $ j\leq 0$.

\vskip3mm

\begin{cor}\label{cor:1}
Let $s \in \mathbb R$. The operators $A^{s/2}$ on 
$\mathcal Z(A)$ and $\mathcal Z'(A)$ are well-defined and continuous. 
\end{cor}

\vskip3mm

We will prove in Proposition \ref{Prop.2.6} below that the semigroup generated by $-A^{\alpha/2}$, $\alpha>0$, is well-defined on test function spaces. This well-definedness will allow to define the semigroup on both
$\mathcal X'(A)$ and $\mathcal Z'(A)$ as corresponding dual operators. 

To define the semigroup generated by the fractional operator $A^{\alpha/2}$ rigorously, we utilize subordination (see \cite[Chapter IX.11, pp. 259-268]{Yosida_1980}). 
For $\sigma > 0 $, $t > 0$, $0 < \alpha < 1$, 
let us introduce the subordinator $F_{t,\alpha}$, 
\begin{equation}\label{0829-1}
F_{t,\alpha} (s) 
:= 
\begin{cases}
\displaystyle 
\frac{1}{2\pi i} \int_{\sigma - i \infty} ^ {\sigma + i \infty} 
 e^{z s - t z^\alpha } \,{\rm d}z , 
 & s \geq 0, 
\\
0, 
& s < 0,
\end{cases}
\end{equation}
where the branch of $z^\alpha$ is taken so that ${\rm Re } \, z^\alpha > 0$ 
provided that ${\rm Re \,} z > 0$. 
It is easily seen by the residue theorem that $F_{t,\alpha}$ is independent of 
$\sigma > 0$, 
by a change of variable 
\[
F_{t,\alpha} (s) = t^{-\frac{1}{\alpha}} F_{1,\alpha} (t^{-\frac{1}{\alpha}} s) ,
\]
and by a careful estimate with choosing $\sigma = \frac{1}{s}$ 
\[
F_{t,\alpha}(s) = o(s^N) 
\quad \text{as } s \to 0 
\quad \text{ for all } N \in \mathbb N . 
\]
As shown in~\cite{Yosida_1980}, $F_{t,\alpha}$ satisfies the following fundamental properties: 
\begin{align} 
& F_{t,\alpha} \geq 0,\nonumber\\  
& \int _{0}^\infty F_{t,\alpha} (s) \,{\rm d}s = 1 ,\label{conv}\\
& F_{t_1 , \alpha} * F_{t_2 , \alpha} (s) =F_{t_1 + t_2 , \alpha} (s).\nonumber
\end{align}
Then the semigoup generated by the fractional power $-A^{\alpha/2}$
is introduced by means of the equality (see \cite[Proposition 1, p. 260]{Yosida_1980})
\begin{equation}\label{0806-2}
e^{-t \lambda ^\alpha} = \int_0^\infty F_{t,\alpha} (s) e^{-s\lambda } \,{\rm d}s , 
\quad \lambda  > 0 .
\end{equation}

This formula allows to handle the fractional case based on the non-fractional one. 

\vskip3mm

\begin{prop}\label{Prop.2.6}
The semigroups $\{ e^{-t A^{\alpha/2}}\} _{ t\geq 0}$ on $\mathcal X(A)$, $\mathcal X'(A)$, 
$\mathcal Z(A)$, and $\mathcal Z'(A)$ are well-defined and uniformly bounded with respect to $t$. 
\end{prop}

\begin{pf} 
We prove the well-definedness on $\mathcal X(A)$. 
Let $\ell_0 \in \mathbb N$ be such that $2\ell_0 > \alpha$. 
Let $f \in \mathcal X(A)$. 
Since $f \in L^2 (\Omega)$, $e^{-tA^{\ell_0}} f$ is initially defined by the 
spectral decomposition 
\[
e^{-tA^{\ell_0}}f = \int _0^\infty e^{-t \lambda ^{\ell_0}} \,{\rm d} E(\lambda) f 
\quad \text{ in } L^2 (\Omega) .
\]
Since $\ell_0 \in \mathbb N$,  $e^{-t \lambda ^{l_0}}$ ($\lambda \geq 0$) can be extended to a function in the Schwartz class 
on the real line, and the boundedness of the spectral multiplier \eqref{0806-1} implies that
\begin{equation}\label{0826-1}
\sup _{t > 0} \| e^{-tA^{\ell_0}} f \|_{L^1 (\Omega)} 
\leq C \| f \|_{L^1(\Omega)}.
\end{equation}
Therefore, 
\[
p_m (e^{-t A^{l_0}}f) \leq C\, p_m(f) 
\quad \text{for all }  m \in \mathbb N, \; t\geq 0,  
\text{ and fixed } \ell_0 \in \mathbb N. 
\]
and $\{ e^{-t A^{l_0}}\} _{ t\geq 0}$ can be regarded as a family of uniformly bounded operators on $\mathcal X(A)$.
In the case of fractional powers, note that 
\[
\lambda ^{\alpha/2} = (\lambda ^{l_0} ) ^{\frac{\alpha}{2 l_0}}, 
\qquad 
e^{-t \lambda^{\alpha /2}} 
= \int_0 ^\infty F_{t,\frac{\alpha}{2 l_0}} (s) e^{-s \lambda^{l_0}} \,{\rm d}s , 
\quad \lambda > 0.
\]
The latter identity leads via (\ref{0806-2}) to an equality for the semigroup on $L^2 (\Omega)$,
\[
\begin{split}
e^{-t A^{\alpha /2}} 
=
& \int_0 ^\infty e^{-t \lambda ^{\alpha /2}}\,{\rm d}E(\lambda) 
= \int_0 ^\infty 
 \Big(  \int_0 ^\infty F_{t,\frac{\alpha}{2 l_0}} (s) e^{-s \lambda^{ l_0}} \,{\rm d}s
 \Big) 
\,{\rm d}E(\lambda) 
\\
= 
& 
 \int_0 ^\infty F_{t,\frac{\alpha}{2 l_0}} (s) e^{-s A^{ l_0}} \,{\rm d}s,
\end{split}
\]
together with the estimate 
\begin{equation}\label{0826-2}
\sup_{t > 0} \| e^{-t A^{\alpha /2}} f \|_{L^1(\Omega)} 
\leq \sup_{t > 0}  \int_0 ^\infty F_{t,\frac{\alpha}{2 l_0}} (s) 
  \|e^{-s A^{ l_0}}f \|_{L^1 (\Omega)} \,{\rm d}s
\leq 
C \| f \|_{L^1 (\Omega)},
\end{equation}
where the constant $C>0$ depends on $\ell_0$. However, we remark that the representation of the semigroup on $L^2$ does not depend on the choice of $\ell_0$ satisfying $2\ell>\alpha_0$.
In view of \eqref{0826-1} and its consequence on $p_m (e^{-t A^{l_0}}f)$ we conclude from \eqref{0826-2} that  
\[
p_m (e^{-t A^{\alpha}}f) \leq C\, p_m(f) 
\quad \text{for all }  m \in \mathbb N . 
\]
Consequently, $\{e^{-t A^{\alpha /2}}\}_{t\geq 0}$ is a bounded operator family on $\mathcal X(A)$ 
for all $\alpha > 0$, since $l_0$ is arbitrary. 
By duality, we also get that it is well-defined on $\mathcal X'(A)$. 

The well-definedness and boundedness on $\mathcal Z(A)$ and $\mathcal Z'(A)$ follow analogously. 

In each case, the semigroup property of $\{e^{-t A^{\alpha /2}}\}_{t\geq 0}$ is based on the convolution property $\eqref{conv}_3$.
\end{pf}

\vskip3mm 

We now fix the definition of the semigroup generated by $A^{\alpha /2}$ on dual spaces. 

\begin{definition}\label{def:2}
(Fractional semigroups) 
Let $\alpha > 0$. 
\begin{enumerate}
\item[(i)] 
The semigroup $\{ e^{-t A^{\alpha/2}}\} _{ t\geq 0}$ on $\mathcal X'(A)$ is defined by 
\[
_{\mathcal X'}\langle  e^{-t A^{\alpha/2}} f , g \rangle _{\mathcal X} 
:= \, _{\mathcal X'}\langle  f , e^{-t A^{\alpha/2}} g \rangle _{\mathcal X} , 
\quad f \in \mathcal X'(A), g \in \mathcal X(A) .
\]
\item[(ii)] 
By analogy, the semigroup $\{ e^{-t A^{\alpha/2}}\} _{ t\geq 0}$ on $\mathcal Z'(A)$ is defined by 
\[
_{\mathcal Z'}\langle  e^{-t A^{\alpha/2}} f , g \rangle _{\mathcal Z} 
:= \, _{\mathcal Z'}\langle  f , e^{-t A^{\alpha/2}} g \rangle _{\mathcal Z} , 
\quad f \in \mathcal Z'(A), g \in \mathcal Z(A). 
\]
\end{enumerate}
\end{definition}

\vskip3mm 

We recall the definition of Besov spaces based on the operator $A$ (see~\cite{IMT-2019}). Here we also need the partition of unity with the function $\psi$, see \eqref{0301-1}.

{\begin{definition}\label{D2.7} 
Let $s \in \mathbb R$ and $1 \leq p,q \leq \infty$. 
\begin{enumerate}
\item[(i)] 
The Besov spaces of non-homogeneous type are defined by 
\[
B^s_{p,q} (A) := 
\left\{ f \in \mathcal X'(A) \, | \, 
   \| f \|_{B^s_{p,q} (A)} < \infty 
\right\},
\]
where 
\[
\| f \|_{B^s_{p,q} (A)} 
:= \| \psi (\sqrt{A}) f \|_{L^p} 
  + \left\| \Big\{ 2^{sj} \| \phi_j (\sqrt{A}) f \|_{L^p (\Omega)} 
            \Big\} _{j \in \mathbb N}
  \right\|_{\ell^q(\mathbb N)} . 
\]

\item[(ii)] 
The Besov spaces of homogeneous type are defined by 
\[
\dot B^s_{p,q} (A) := 
\left\{ f \in \mathcal Z'(A) \, | \, 
   \| f \|_{\dot B^s_{p,q} (A)} < \infty 
\right\},
\]
where 
\[
\| f \|_{\dot B^s_{p,q} (A)} 
:= \left\| \Big\{ 2^{sj} \| \phi_j (\sqrt{A}) f \|_{L^p (\Omega)} 
            \Big\} _{j \in \mathbb Z}
  \right\|_{\ell^q(\mathbb Z)} . 
\]
\end{enumerate}
\end{definition}

We prepare two lemmata below for the proof of the main theorems in the following sections. 

\begin{lem}\label{lem:0912-1}
Let $s \in \mathbb R$ and $1 \leq p \leq \infty$. Then
\[
B^s_{p,1} (A) \hookrightarrow H^s_p (A) 
\hookrightarrow  B^s_{p,\infty}(A),
\qquad 
\dot B^s_{p,1} (A) \hookrightarrow \dot H^s_p (A) 
\hookrightarrow \dot B^s_{p,\infty}(A) .
\]
\end{lem}

\vskip3mm 

\noindent
{\bf Proof.}
We only consider the homogeneous case, since the non-homogeneous case is similar. 
Let $f \in \dot B^s_{p,1} (A)$. Then we get by the resolution of the identity in $\mathcal Z'(A)$ \BLACK  (Lemma~\ref{lem:5}), 
the triangle inequality and Lemma~\ref{lem:1} that
\[
\| f \|_{\dot H^s_p (A)} 
\leq \sum_{ j \in \mathbb Z} \| \phi_j (\sqrt{A}) f \|_{\dot H^s_p (A)} 
\leq C\sum_{ j \in \mathbb Z} 2^{sj} \| \phi_j (\sqrt{A}) f \|_{L^p}  
= \| f \|_{\dot B^s_{p,1}(A)} ,
\]
which proves the first embedding $B^s_{p,1} (A) \hookrightarrow H^s_p (A)$. 
As for the second one, 
\[
\| f \|_{\dot B^s_{p,\infty}(A)} 
\leq C \sup_{j \in \mathbb Z} \| \phi_j (\sqrt{A}) f \|_{\dot H^s_p (A)} 
\leq C \| f \|_{\dot H^s_p (A)},
\]
where in the first step we used \eqref{0802-1} to see that $2^{sj}\|\phi_j(\sqrt{A}) f\|_{L^p} \leq c\|A^{s/2}\phi_j(\sqrt{A})f\|_{L^p}$. 
Hence, we conclude the second embedding 
$\dot H^s_p (A) \hookrightarrow \dot B^s_{p,\infty}(A)$. 
\hfill $\Box$ 

\vskip3mm

\begin{lem}\label{lem:0912-2} 
Let $s \in \mathbb R$, $1 \leq p \leq \infty$ and $f \in H^s_p(A)$. Then
\begin{equation}\label{0925-1}
\| f \| _{H^{s}_{p}(A)}
\simeq 
\| \psi (\sqrt{A})f \| _{L^{p}}
+ \Big\| \sum _{ j=1}^\infty \phi_j (\sqrt{A}) f \Big\|_{\dot H^{s}_{p}(A)}.
\end{equation}
\end{lem}

\vskip1mm 

\noindent
{\bf Proof.}
Let $f \in H^s_{p} (A)$. By Lemma~\ref{lem:5}, we decompose $f$ into 
\[
f = \psi(\sqrt{A}) f + \sum_{j=1}^\infty \phi_j (\sqrt{A})  f 
=: f_L + f_H
\quad \text{in } \mathcal X'(A) .  
\]
It follows from Lemma~\ref{lem:1} that 
\[
\| f_L \|_{H^s_p (A)} 
= \| (1+A)^{s/2}\psi (2^{-2}\sqrt{A}) f_L \|_{L^p}
\leq C \| f_L \|_{L^p}
\]
since $\psi (\sqrt{A}) = \psi (2^{-2}\sqrt{A}) \psi (\sqrt{A})$ 
and $(1+\lambda)^{s/2}\psi (2^{-2}\sqrt{\lambda})$ is a smooth function with compact support. 
As for the high spectral component, 
by the triangle inequality, 
\[
\begin{split}
\| f_H \|_{H^s_p} 
\leq 
\| A^{s/2} f_H \|_{L^p} + \| ((1+A)^{s/2} - A^{s/2} )f_H \|_{L^p}.
\end{split}
\] 
We write by the resolution of the identity, see Lemma~\ref{lem:5}, and the spectrum of $f_H$ 
away from the origin that 
\[
\begin{split}
\Big( (1+A)^{s/2} - A^{s/2} \Big)f_H 
=& \sum_{j=0}^\infty \phi_j (\sqrt{A}) \Big( (1+A)^{s/2} - A^{s/2} \Big) 
A^{-s/2} ( A^{s/2} f_H ) 
\quad \text{ in } \mathcal X'(A). 
\\
=:& \sum_{j=0}^\infty \widetilde \phi_j (\sqrt{A}) (A^{s/2}f_H). 
\end{split}
\]
By Lemma~\ref{lem:1} and the definition $\widetilde  \phi_{0,j}(\lambda) := \widetilde \phi_j(2^{j}\lambda) $ we get that 
\[
\| \widetilde \phi_j (\sqrt{A}) (A^{s/2}f_H) \|_{L^p} 
\leq C \| (1+|\cdot|^2)^{\frac{\gamma}{2}} \widetilde \phi_{0,j} \|_{H^\beta} 
\| A^{s/2} f_H\|_{L^p}, 
\]
where $\beta > (d+1)/2$ ($\beta\in\mathbb N$) and $\gamma > (3d+2)/4$. 
To estimate the $H^\beta$ norm we exploit the fundamental theorem of 
calculus and get
\[
\widetilde  \phi_{0,j} (\lambda)
=  \phi_{0} (\lambda) \int_0^1 \partial _\theta ( \theta +2^{2j}\lambda^2)^{s/2}   d\theta 
  \, (2^j  \lambda)^{-s}
= \frac{s2^{-2j}}{2 \lambda^{s/2}} \phi_0(\lambda) \int_0^1 ( 2^{-2j}\theta +\lambda^2)^{s/2 -1} d\theta. 
\]
Since the compact support of $\phi_0$ is away from the origin, it is easy to see that 
\[
\| (1+|\cdot|^2)^{\frac{\gamma}{2}}  \widetilde  \phi_{0,j}\BLACK  \|_{H^\beta}
\leq C 2^{-2j} .
\]
Therefore 
\[
\sum_{j=0}^\infty\| \widetilde \phi_j (\sqrt{A}) (A^{s/2}f_H) \|_{L^p} 
\leq C \sum_{j=0}^\infty 2^{-2j} 
\| A^{s/2} f_H\|_{L^p} 
= C \| f_H \|_{\dot H^s_p(A)} . 
\]
which verifies that the left hand side of \eqref{0925-1} is controlled by the right hand side. 
\BLACK 

For the converse inequality, we apply a similar argument. In fact, it follows that 
\begin{align*}
\| f_L \|_{L^p} 
& =\| (1+A)^{-s/2}\psi (\sqrt{A}) (1+A)^{s/2}f \|_{L^p}  
\\
& 
\leq C \| (1+A)^{s/2} f \|_{L^p} 
= C\|f\|_{H^s_p(A)} ,
\end{align*}
and 
\begin{align*}
\| f_H \|_{\dot H^s_p (A)} 
& 
= \| (1+A) ^{s/2} f_H \|_{L^p} 
 + \| ((1+A)^{s/2} - A^{s/2} ) f_H \|_{L^p}
\\
& \leq 
C \| f \|_{H^s_{p}(A)} .
\end{align*}
\hfill $\Box$

Finally, we prove that the semigroup and the fractional Laplacian commute. 

\begin{lem}\label{lem:0301-2}
Let $\alpha > 0$. 
\begin{enumerate}
\item[(i)] 
For $s > 0$, $A^{s/2}$ and $e^{-t A^{\alpha /2}}$ commute on $\mathcal X'(A)$. 
\item[(ii)] 
For $s \in \mathbb R$, 
$(1+A)^{s/2}$ and $e^{-t A^{\alpha /2}}$ commute on $\mathcal X'(A)$. 
\item[(iii)]
For $s \in \mathbb R$, 
$A^{s/2}$ and $e^{-t A^{\alpha /2}}$ commute on $\mathcal Z'(A)$. 
\end{enumerate}
\end{lem}

\noindent 
{\bf Proof. } We only prove (iii), since (i) and (ii) can be proved by a similar argument. 
Let us start by proving that for $f \in \mathcal Z(A)$
\[
A^{s/2} e^{-t A^{\alpha /2}} f =  e^{-t A^{\alpha /2}}A^{s/2} f 
\text{ in } \mathcal Z(A) .
\] 
Indeed, since $f \in \mathcal Z(A) \subset L^2$ and $\lambda^{s/2}\cdot e^{-t \lambda^{\alpha /2}} =e^{-t \lambda^{\alpha /2}}\cdot \lambda^{s/2}$  as scalar valued functions, we get the equality by the spectral theorem
\[
\begin{split}
A^{s/2} e^{-t A^{\alpha /2}} f 
= & \int _0^\infty  \lambda^{s/2} e^{-t \lambda^{\alpha /2}}\,{\rm d} E(\lambda) f
=  \int _0^\infty  e^{-t \lambda^{\alpha /2}}\lambda^{s/2}\,{\rm d} E(\lambda) f
\\
= &\, e^{-t A^{\alpha /2}}A^{s/2} f \;\text{ in } L^2.
\end{split}
\]
Therefore 
$A^{s/2} e^{-t A^{\alpha /2}} f (x)$ equals
$e^{-t A^{\alpha /2}} A^{s/2} f (x)$ for almost every $x$. 
We then conclude that $A^{s/2} e^{-t A^{\alpha /2}} f =e^{-t A^{\alpha /2}} A^{s/2}  f $ 
in $\mathcal Z(A)$. 
The commutativity in $\mathcal Z'(A)$ 
follows from  commutativity in $\mathcal Z(A)$ and a duality argument. 
\hfill $\Box$

\section{Proof of Theorem~\ref{thm:1}}

\noindent 
{\bf Proof of (i). }  
We start by the completeness of $\dot H^s_p (A)$. 
Let $\{ f_n \}_{n=1}^\infty$ be a Cauchy sequence in $\dot H^s_p (A)$. 
By definition, 
$\{ A^{s/2} f_n \}_{n=1}^\infty$ is a Cauchy sequence in $L^p (\Omega)$ 
and $A^{s/2} f_n$ converges to an element $g \in L^p (\Omega)$ 
as $n \to \infty $ due to the completeness of $L^p (\Omega)$. 
We have $g, A^{-s/2} g  \in \mathcal Z'(A)$, 
since $L^p (\Omega) \hookrightarrow \mathcal Z'(A)$ (see Lemma~\ref{lem:0301-3}) 
and $A^{-s/2}$ is well-defined on $\mathcal Z'(A)$ (see Corollary~\ref{cor:1}). 
We also see that $f_n \to A^{-s/2} g$ in $\dot H^s_p (A)$, 
which proves the completeness of $\dot H^s_p(A)$. 
It is easy to see that $\| \cdot \|_{\dot H^s_p (A)}$ is a norm, 
and hence, $\dot H^s_p (A) $ is a Banach space. 

We turn to prove the embedding 
$\mathcal Z(A) \hookrightarrow \dot H^s_p (A) \hookrightarrow \mathcal Z'(A)$. 
Let $m \in \mathbb N$ be such that $m > |s + d(1-1/p)|$. 
For $f \in \mathcal Z(A)$, we have from Lemma \ref {lem:5} and 
\eqref{0802-1} that 
\[
\| f \|_{\dot H^s_p (A)} 
\leq C \sum _{j \in \mathbb Z} 2^{s j} \| \phi_j (\sqrt{A}) f \|_{L^p (\Omega)}
\leq C \sum _{j \in \mathbb Z} 2^{s j + d(1-\frac{1}{p})j} 
       \| \phi_j (\sqrt{A}) f \|_{L^1 (\Omega)}
\leq C q_{m} (f) ,
\]
which proves $\mathcal Z(A) \hookrightarrow \dot H^s_p (A)$.
 
For the second embedding $\dot H^s_p (A) \hookrightarrow \mathcal Z'(A)$, choose $m' > |-s + d/p|$. 
For $f \in \dot H^s_p (A)$, 
it follows from \eqref{0802-1} and the previous inequality for any $g\in \mathcal Z(A) \subset\dot H^{-s}_{p'}(A)$ that 
\[
\big| _{\mathcal Z'} \langle f ,g \rangle _{\mathcal Z} \big| 
= \big| _{\mathcal Z'} \langle A^{s/2} f , A^{-s/2}g  \rangle _{\mathcal Z} \big| 
\leq \| f \|_{\dot H^s_p (A)} \| g \|_{\dot H^{-s}_{p'} (A)}
\leq C \| f \|_{\dot H^s_p (A)} q_{m'} (g).
\]

The proof for non-homogeneous spaces $\mathcal X(A)$, $\mathcal X'(A)$, and
$H^s_p(A)$ is similar. 
\hfill $\Box$ 

\vskip3mm

\noindent 
{\bf Proof of (ii). } 
We prove $\dot H^{-s}_{p'} (A) \subset (\dot H^s_{p} (A))' $. 
For any $f \in \dot H^{-s}_p (A)$, we define the linear operator 
$T_f : \dot H^s_p (A) \to \mathbb C$ by 
\[
T_f (g):= 
\int_{\Omega} \big(A^{-s/2}f \big ) \,  \overline{A^{s/2} g} \,{\rm d}x, 
\quad g \in \dot H^s_p(A) .
\]
Then $T_f$ is continuous and the H\"older inequality gives 
\[
|T_f (g)| \leq \| f \|_{\dot H^{-s}_{p'}(A)} \| g \|_{\dot H^s_p (A)}.
\]

Conversely, 
for any $f \in (\dot H^s_p (A)) '$, we define the linear operator 
$T_f'\BLACK  : L^p (\Omega) \to \mathbb C$ by 
\[
T_f'(G) := f(A^{-s/2} G) , \quad G \in L^p (\Omega). 
\]
The continuity of $T_f'$ follows from 
\[
|T_f'(G)| 
= |f(A^{-s/2} G)| 
\leq \| f \|_{(\dot H^s_p (A)) '} \| A^{-s/2} G \|_{\dot H^s_{p}(A)}
= \| f \|_{(\dot H^s_p (A)) '} \| G \|_{L^p} .
\]
Since $(L^p(\Omega))' = L^{p'} (\Omega)$, 
there exists $F \in L^{p'}(\Omega)$ such that $\|F\|_{L^{p'}} = \|T_f'\|_{(L^p)'}\BLACK $ and 
\[
T_f' (G) = \int_{\Omega} F(x) \, \overline{G(x)} \,{\rm d}x, 
\quad G \in L^p (\Omega) .
\]
Here for any $g \in \dot H^s_p (A)$, 
\[
f(g) = T_f'(A^{s/2} g) = \int _{\Omega} F(x) \, \overline{A^{s/2} g (x)} \,{\rm d}x.
\]
In particular, by the definition of $A^{s/2}$ on $\mathcal Z'(A)$, 
\[
f(g) 
= \int _{\Omega} F(x) \, \overline{A^{s/2} g (x)} \,{\rm d}x 
= _{\mathcal Z'} \langle F , A^{s/2} g \rangle _{\mathcal Z}
= _{\mathcal Z'} \langle A^{s/2} F , g \rangle _{\mathcal Z}, 
\quad g \in \mathcal Z(A). 
\]
Then $\widetilde F := A^{s/2} F \in \dot H^{-s}_{p'} (A)$
and $\|\widetilde F \|_{\dot H^{-s}_{p'}(A)} 
= \| F \|_{L^{p'}} = \| T_f \|_{(L^p)'} \leq \| f \|_{(\dot H^s_p (A))'}$,
which proves $(\dot H^s_{p} (A))' \subset \dot H^{-s}_{p'} (A) $. 

The duality for the non-homogeneous spaces follows analogously. 
\hfill $\Box$ 

\vskip3mm

\noindent 
{\bf Proof of (iii). } 
The lift property follows from the well-definedness of 
$(1+A)^{s_0 /2}$, $A^{s_0 /2}$ on  $\mathcal X'(A)$, $\mathcal Z'(A)$, respectively, 
and the definition of the norm of $\dot H^s_p (A)$. The same arguments apply to the non-homogeneous spaces.
\hfill $\Box$ 

\vskip3mm

\noindent 
{\bf Proof of (iv). } 
It suffices to prove the case when $s = 0$ by the lift property. 
We apply the argument for the proof of the equivalence of decay estimates of 
semigroups and Sobolev inequalities, {\em cf.} the proof of Theorem~2.4.2 in the monograph 
by Davies~\cite{Davies_1990}. Putting 
\[
\gamma := d \Big( \frac{1}{r} -\frac{1}{p} \Big), 
\]
we write $f\in\dot H^\gamma_r(A)$ in the form
\[
f = A^{-\gamma/2} A^{\gamma/2} f 
= c \int_0 ^\infty t^{-1+\frac{\gamma}{2}}e^{-tA} (A^{\gamma/2} f) \,{\rm d}t 
\quad \text{ in } \mathcal Z'(A)
\]
with the constant $c=\Gamma(\gamma/2) > 0$ and decompose it into 
$f = f_1 + f_2$ where 
\[
f_1 := c \int_0 ^T t^{-1+\frac{\gamma}{2}}e^{-tA} (A^{\gamma/2} f) \,{\rm d}t , \quad 
f_2 := c \int_T ^\infty t^{-1+\frac{\gamma}{2}}e^{-tA} (A^{\gamma/2} f) \,{\rm d}t ,
\]
with $T >0$ to be fixed later. Finally, we recall that the distribution function 
$
\mu _f (\lambda) = | \{ x \,| \, |f (x)| > \lambda \} |
$ 
satisfies the elementary inequality
\[
\mu_{f} (\lambda) 
\leq \mu_{f_1} \Big(\frac{\lambda}{2}\Big) + \mu _{f_2} \Big( \frac{\lambda}{2} \Big) .
\]
Now we use the Gaussian estimate of the semigroup $\{e^{-tA}\}_{t\geq 0}$, see \cite[Example~2.1.8]{Davies_1990}, to get by H\"older's inequality the estimate
\[
\| f_2 \|_{L^\infty} 
\leq C\int_T ^\infty t^{-1+\frac{\gamma}{2}}t^{-\frac{d}{2r}} \| A^{\gamma/2} f \|_{L^r} \,{\rm d}t
\leq C T ^{\frac{\gamma}{2} - \frac{d}{2r}} \| A^{\gamma/2}f \|_{L^r}
= C T ^{-\frac{d}{2p}} \| A^{\gamma/2}f \|_{L^r}. 
\]
We choose $T > 0$ such that 
\[
\frac{\lambda}{2} = C T ^{-\frac{d}{2p}} \| A^{\gamma/2}f \|_{L^r}, 
\]
and get $\mu_{f_2} (\lambda /2) = 0$. 
This and the boundedness of $e^{-tA}$ imply 
\[
\mu _f (\lambda) 
\leq \mu_{f_1} \Big( \frac{\lambda}{2} \Big) 
\leq \frac{2^r}{\lambda ^r} \| f_1 \|_{L^r} ^r
\leq \frac{2^r}{\lambda ^r}  
     \Big( C T^{\frac{\gamma}{2}} \| A^{\gamma/2} f\|_{L^r}\Big) ^r 
= C \lambda ^{-p} \| A^{\gamma /2} f \|_{L^r} ^p, 
\]
which proves that the the operator $A^{-\gamma /2}$ is of weak type $(r,p)$.
Then the Marcinkiewicz interpolation theorem yields
the boundedness of $A^{-\gamma /2} $ from $L^r (\Omega)$ to 
$L^p (\Omega)$ provided that $1 < r < p < \infty$. 
Therefore, we obtain $\dot H^{\gamma}_r (A) \hookrightarrow L^p (\Omega)$. 

The embedding for the non-homogeneous case follows from a similar argument 
by replacing $A$ with $1 + A$. 
\hfill $\Box$ 

\vskip3mm  

\noindent 
{\bf Proof of (v). }
Put 
\[
\dot Y^s_p := \Big\{ f \in \mathcal X'(A) \, \Big| \, 
\| f \|_{\dot H^s_p (A)} < \infty , \,\,
 f = \sum _{j \in \mathbb Z} \phi_j (\sqrt{A}) f \text{ in } \mathcal X'(A)
\Big\}, 
\]
\[
\| f \|_{\dot Y^s_p} :=
\Big\| \sum _{j \in \mathbb Z} \phi_j (\sqrt{A}) f \Big\|_{\dot H^s_p (A)}.
\]
It follows from  $\mathcal X'(A) \hookrightarrow \mathcal Z'(A)$ that 
$\dot Y^s_p \hookrightarrow \dot H^s_p(A)$. For the converse, we consider the map 
\[
T: \dot H^s_p (A) \ni f \mapsto 
\sum _{j \in \mathbb Z} \phi_j (\sqrt{A}) f \in \dot Y^s_p .
\]
This map is well-defined. In fact, for $f \in \dot H^s_p (A)$, it follows from \eqref{0802-1} that 
\[
\sum _{j \leq 0} \| \phi_j (\sqrt{A}) f\|_{L^\infty} 
\leq C \sum _{j \leq 0}  2^{\frac{d}{p}j} \| \phi_j (\sqrt{A}) f \|_{L^p (\Omega)} 
\leq C \Big( \sum _{j \leq 0}  2^{(\frac{d}{p} - s) j} \Big) \| f \|_{\dot H^s_p (A)} 
< \infty ,
\]
which proves that the low spectral part of $f$ belongs to 
$L^\infty (\Omega) \subset \mathcal X'(A)$. For the high spectral part,  
by the definition of $\mathcal X'(A)$ it is easy to see that 
\[
\sum _{ j \geq 1} \phi_j (\sqrt{A}) f \in \mathcal X'(A) ,
\]
since 
\[  
\Big\|\sum_{j\geq 1} \phi_j(\sqrt{A})f\Big\|_{L^p} \leq C\sum_{j\geq 1} 2^{-sj} \|A^{s/2}f\|_{L^p} = c\|f\|_{\dot H^s_p(A)}. 
\]
%
%
Thus, $f \in \mathcal X'(A)$. 
By definition, $T$ is an isometry. 
Therefore, we obtain $\dot H^s_p (A) \hookrightarrow \dot Y^s_p$. 
\hfill $\Box$

\section{Proof of Theorem~\ref{thm:2}}

\noindent
We start by proving the inequality in Besov spaces, namely. 
\begin{equation}\label{0802-3}
\| f \|_{\dot B^s_{p,1}(A)} 
\leq C \| f \|_{\dot B^0_{r,\infty}(A)} ^\theta \| f \|_{\dot B^{s_0}_{r_0,\infty}(A)} ^{1-\theta} . 
\end{equation}
Once it is proved, the inequality for the Sobolev spaces in Theorem~\ref{thm:2} 
is an immediate consequence of \eqref{0802-3} and the embedding 
$\dot B^s_{p,1} (A) \hookrightarrow \dot H^s_p (A) \hookrightarrow \dot B^s_{p,\infty} (A)$ 
for all $s \in \mathbb R$ and $1 \leq p \leq \infty$, see
Lemma~\ref{lem:0912-1}.

Therefore, let us show inequality \eqref{0802-3}. Hereafter, we write $f_j = \phi_j (\sqrt{A}) f$ for simplicity.

\vskip3mm

\noindent 
\underline{The case when $r ,r_0 \leq p $}. 
We notice that 
$s > 0$ and $r \leq p$ imply $-d/r < s-d/p $, so that  $s-d/p$ as a convex combination in \eqref{0912-1} 
satisfies   
\begin{equation}\label{0926-1}
s- \frac{d}{p} \in \Big( - \frac{d}{r} , s_0 - \frac{d}{r_0} \Big) .
\end{equation}
For $N \in \mathbb Z$ we split the infinite series in the definition of the norm of Besov spaces into two series, one for $j \leq N$ and another for $j > N$. Applying \eqref{0802-1} 
we get that
\[
\begin{split}
\| f \|_{\dot B^s_{p,1}(A)} 
\leq 
& C \sum _{ j \leq N} 2^{sj + d (\frac{1}{r}- \frac{1}{p})j}  
   \| f_j \|_{L^{r}} 
+ C\sum _{ j > N} 2^{sj + d(\frac{1}{r_0}-\frac{1}{p})j} 
  \| f_j \|_{L^{r_0}}
\\
\leq 
& C 2^{sN + d (\frac{1}{r}- \frac{1}{p})N}  
   \| f \|_{\dot B^0_{r,\infty}(A)} 
+  C 2^{sN + d(\frac{1}{r_0}-\frac{1}{p})N - s_0N} 
\| f \|_{\dot B^{s_0}_{r_0,\infty}(A)},
\end{split}
\]
since $s+d(1/r - 1/p) > 0, s + d(1/r_0 - 1/p) -s_0 < 0$. By choosing $N$ such that 
\[
 2^{sN + d (\frac{1}{r}- \frac{1}{p})N}  
   \| f \|_{\dot B^0_{r,\infty}(A)} 
\simeq  2^{sN + d(\frac{1}{r_0}-\frac{1}{p})N - s_0N} 
\| f \|_{\dot B^{s_0}_{r_0,\infty}(A)},
\]
we obtain \eqref{0802-3} with the same $\theta$ as in Theorem~\ref{thm:2} 
under the condition $r_,r_0 \leq p$. 

\vskip3mm 

For the next two cases, let $\mu $ satisfy $1/ p = \mu/r + (1-\mu)/ r_0$.  

\vskip3mm

\noindent 
\underline{The case when $r \leq p < r_0$}. 
As in the previous case, we have \eqref{0926-1}, 
and, due to \eqref{0912-1} and $s < (1-\theta) s_0$, 
\[
 \frac{d}{p} - \frac{d}{r_0} 
  - \theta \Big( \frac{d}{r} - \frac{d}{r_0} \Big) 
  =  s- (1 - \theta) s_0 < 0   
\]
which implies $1/p < \theta /r + (1-\theta) / r_0$, and 
$\theta > \mu$. 
Hence we have $s +d/r - d/p >0$ and $s - (1- \mu )s_0 < s- (1-\theta ) s_0 < 0$.
For $N \in \mathbb Z$ we split the sum into two parts, $j \leq N$ and $j > N$, apply \eqref{0802-1} 
and the interpolation inequality $\|f_j\|_{L^p} \leq \|f_j\|_{L^r}^\mu \|f_j\|_{L^{r_0}}^{1-\mu}$ to get that
\begin{align}\label{4.3a}
\| f \|_{\dot B^s_{p,1}(A)} 
\leq
& C\sum _{ j \leq N} 2^{sj + d (\frac{1}{r}- \frac{1}{p})j}  
   \| f_j \|_{L^{r}} 
+ C\sum_{j > N} 2^{sj} \|f_j\|_{L^{r}} ^{\mu} \|f_j\|_{L^{r_0}} ^{1-\mu}\nonumber
\\
\leq 
& C 2^{sN + d (\frac{1}{r}- \frac{1}{p})N}  
   \| f \|_{\dot B^0_{r,\infty}(A)} 
+  C 2^{sN - (1-\mu)s_0 N} 
  \| f \|_{\dot B^0_{r,\infty}(A)} ^{\mu} \| f \|_{\dot B^{s_0}_{r_0,\infty}(A)}^{1-\mu} .
\end{align}
By choosing $N$ such that 
\[
2^{sN + d (\frac{1}{r}- \frac{1}{p})N}  
   \| f \|_{\dot B^0_{r,\infty}(A)} 
\simeq  2^{sN - (1-\mu)s_0 N} 
  \| f \|_{\dot B^0_{r,\infty}(A)} ^{\mu} \| f \|_{\dot B^{s_0}_{r_0,\infty}(A)}^{1-\mu}, 
\]
we obtain \eqref{0802-3} with the same $\theta$ as in Theorem~\ref{thm:2}. 

\vskip3mm

\noindent 
\underline{The case when $r_0 \leq p < r$}. 
We start by the case when 
$-d/r < s_0 - d/r_0$.
As in \eqref{4.3a} we obtain that
\begin{equation}\label{0927-2}
\begin{split}
\| f \|_{\dot B^s_{p,1}(A)} 
\leq 
& C\sum _{ j \leq N} 2^{sj}  
   \| f_j \|_{L^r} ^{\mu } \| f_j \|_{L^{r_0}} ^{1-\mu}
+ C\sum _{ j > N} 
 2^{sj + d (\frac{1}{r_0}- \frac{1}{p})j}  
   \| f_j \|_{L^{r_0}} 
\\
\leq 
& C 2^{sN - (1-\mu) s_0 N}  
   \| f \|_{\dot B^0_{r,\infty}(A)} ^{\mu} \| f \|_{\dot B^{s_0}_{r_0,\infty}(A)}^{1-\mu}
+  C 2^{( s-s_0 + d(\frac{1}{r_0}-\frac{1}{p}) ) N} 
  \| f \|_{\dot B^{s_0}_{r_0,\infty}(A)}
\end{split}
\end{equation}
provided that $s - (1- \mu) s_0 > 0$ and $s-s_0 - d/p + d/r_0  < 0$. 
Let us check these two conditions on indices. 
Since $\theta < 1- s/ s_0$ and $s-d/p$ is a convex combination of $-d/r$ and $s_0-d/r_0$ with parameter $\theta$, see \eqref{0912-1},  we notice that 
\[
-\frac{d}{r} < s_0 - \frac{d}{r_0} 
\quad \text{ implies } \quad 
s-\frac{d}{p} >
\Big( 1- \frac{s}{s_0} \Big) \Big( - \frac{d}{r} \Big) + \frac{s}{s_0} \Big( s_0 - \frac{d}{r_0} \Big), 
\]
which is equivalent to 
\[
\frac{1}{p} < \frac{1-\frac{s}{s_0}}{r} + \frac{\frac{s}{s_0}}{r_0}, 
\quad \text{ i.e. } \quad 
s - (1-\mu) s_0 > 0 . 
\]
The second condition $s-d/p - s_0+d/r_0  <$ is satisfied, since $-d/r < s_0 - d/r_0$ 
and $s- d/p$ is given by a convex combination, see \eqref{0912-1}. 
By choosing 
$N$ such that 
\begin{equation}\label{0927-1}
2^{sN - (1-\mu)s_0 N} \| f \|_{\dot B^0_{r,\infty}(A)} ^\mu 
   \| f \|_{\dot B^{s_0}_{r_0,\infty}(A)} ^{1-\mu} 
 \simeq 2^{(s-s_0)N + d(\frac{1}{r_0}-\frac{1}{p})N } 
   \| f \|_{\dot B^{s_0}_{r_0,\infty}(A)} ,
\end{equation}
we obtain \eqref{0802-3}. 

Conversely, if $-d/r > s_0 - d/r_0$, then this implies 
\[
s - s_0 + \frac{d}{r_0} - \frac{d}{p}> 0 
\quad \text{and} \quad 
s- (1-\mu)s_0 < 0.
\]
We estimate, exchanging the inequalities for the low and the high spectral parts 
in \eqref{0927-2}, 
\[
\begin{split}
\| f \|_{\dot B^s_{p,1}(A)} 
\leq 
& 
\sum _{ j \leq N} 
 2^{sj + d (\frac{1}{r_0}- \frac{1}{p})j}  
   \| f_j \|_{L^{r_0}} 
+ 
\sum _{ j > N} 2^{sj}  
   \| f_j \|_{L^r} ^{\mu } \| f_j \|_{L^{r_0}} ^{1-\mu}
\\
\leq 
&  C 2^{( s - s_0 + d(\frac{1}{r_0}-\frac{1}{p}) ) N} 
  \| f \|_{\dot B^{s_0}_{r_0,\infty}(A)}
  + 
   C 2^{sN - (1-\mu) s_0 N}  
   \| f \|_{\dot B^0_{r,\infty}(A)} ^{\mu} \| f \|_{\dot B^{s_0}_{r_0,\infty}(A)}^{1-\mu} .
\end{split}
\]
Choosing $N$ satisfying \eqref{0927-1}, we obtain \eqref{0802-3}. 
\hfill $\Box$ 
%
%

\section{Proof of Theorem~\ref{thm:3}}

\noindent 
{\bf Proof of Theorem~\ref{thm:3} (i), homogeneous case. } 
For each $\alpha > 0$, we know from (\ref{0826-2}) the uniform boundedness of $\{ e^{-tA^{\alpha /2}} \}_{t\geq 0}$ 
in $L^1 (\Omega)$ with respect to $t \geq 0$, which proves the boundedness 
in $L^\infty (\Omega)$ by a duality argument. Then interpolation yields the 
boundedness in $L^p (\Omega)$ for $1 \leq p \leq \infty$. 
Therefore, due to the commutativity of $e^{-tA^{\alpha /2}}$ and $A^{s/2}$, see Lemma \ref{lem:0301-2}, 
\[
\sup_{t > 0} \| e^{-tA^{\alpha /2}} \|_{\dot H^s_p (A) \to \dot H^s_p (A)} 
<\infty 
\quad \text{ for each } \alpha > 0 . 
\]

If the domain is bounded, then the infimum of the spectrum of 
$A$ is strictly positive, which implies that $j_0 \in \mathbb Z$ exists such that 
\[
f =\big( 1 - \psi (2^{-j_0} \sqrt{A} \big) f, 
\quad f \in \dot H^s_p (A) , 
\]
where $\psi$ satisfies \eqref{0301-1}. 
We apply the boundedness of the spectral multiplier \eqref{0806-1} 
to the operator 
$ e^{-t A^{\alpha /2}} \big( 1 - \psi (2^{-j_0} \sqrt{A} \big)$ 
with the scaling of $2^{-j_0}$, and it leads to the exponential decay 
by the following inequality. 
\[
 \| (1+ \lambda^2)^{\beta} 
          e^{-t (2^{j_0} \sqrt{\lambda} ) ^{\alpha}} \big( 1-\psi ( \sqrt{\lambda})\big) \|_{H^\gamma (\lambda \in \mathbb R)} 
\leq C e^{-\kappa t} ,
\]
where $\kappa$ is a positive constant depending on $j_0$ 
and $\beta , \gamma > 1$. 
\hfill $\Box$ 

\vskip2mm

\noindent 
{\bf Proof of Theorem~\ref{thm:3} (ii), homogeneous case. } 
We apply this kind of estimates in Besov spaces 
(see \cite[Theorem~1.1]{Iwa-2018}) and the embedding 
$\dot B^s_{p,1} (A) \hookrightarrow \dot H^s_p(A) \hookrightarrow 
\dot B^s_{p,\infty} (A)$, which imply 
\[
\begin{split}
\big\| e^{-tA^{ \frac{\alpha }{2} }} f \big\|_{\dot H^{s_2}_{p_2} (A)} 
\leq 
& 
C \big\| e^{-tA^{ \frac{\alpha }{2} }} f \big\|_{\dot B^{s_2}_{p_2,1} (A)} 
\leq C 
t^{- \frac{s_2 - s_1}{\alpha}-\frac{d}{\alpha} (\frac{1}{p_1} - \frac{1}{p_2})}  \| f \|_{ \dot B^{s_1}_{p_1,\infty} (A) } 
\\
\leq 
& 
C t^{- \frac{s_2 - s_1}{\alpha}-\frac{d}{\alpha} (\frac{1}{p_1} - \frac{1}{p_2}) }   \| f \|_{ \dot H^{s_1}_{p_1} (A) } .
\end{split}
\]
\hfill $\Box$

\vskip2mm 

\noindent 
{\bf Proof of Theorem~\ref{thm:3} (iii), homogeneous case. } 
%
%
%
It is sufficient to prove the case when $s = 0$ by the lift property 
and the commutativity of $A^{s/2}$ and $e^{-tA^{\alpha/2}}$. 
 
We start by the case when $ 1 \leq p < \infty$. Let $f \in L^p (\Omega)$. 
The ideas are the denseness of $C_0 ^\infty (\Omega)$ in $L^p (\Omega)$ and the continuity of $e^{-tA^{\ell_0}}$ at $t=0$ for the power $\ell\in\mathbb N$ together with subordination. 
In fact, 
since $C_0 ^\infty (\Omega)$ is dense in $L^p (\Omega)$, 
choose a sequence $\{ f_n \}_{n=1}^\infty \subset C_0^\infty (\Omega)$ 
such that $f_n$ converges to $f$ in $L^p (\Omega)$ as $n \to \infty$. 
Then, by the uniform boundedness of $e^{-t A^{\alpha /2}}$ in $L^p(\Omega)$ with respect to $t>0$, 
it suffices to show 
the convergence for $f_n$. We call $f_n \in C_0^\infty(\Omega)$ 
again $f$ and prove that 
\[
\lim _{ t\to 0} \| e^{-t A^{\alpha /2} } f - f \|_{L^p} = 0  ,
\quad f\in C_0 ^\infty (\Omega). 
\]
We write the decomposition of the identity in the form 
\[
\psi (2^{-N} \sqrt{A}) + \sum _{j \geq N}  \phi_j (\sqrt{A}), 
\]
where $\psi, \phi_j$ satisfy \eqref{0301-1}, to analyze low and high spectrum parts separately.
Since $f \in H^s_p(A)$ for all $s > 0$, it follows that 
\begin{align*}
& \sup _{ t > 0} \Big\|  \sum _{j \geq N}  \phi_j (\sqrt{A}) 
   \Big(e^{-t A^{\alpha /2} } f - f  \Big) \Big\|_{L^p}  
\leq 
 C \Big\|  \sum _{j \geq N}  \phi_j (\sqrt{A})  f \Big\|_{L^p}\\
\RED 
&\quad \leq C \sum_{ j \geq N } 2^{-sj} \| A^{s/2} f \|_{L^p}
\leq C 2^{-sN} \| f \|_{H^s_p (A)} \to 0
\end{align*}
as $N \to \infty$. For the low spectrum part, we need to show 
\begin{equation}\label{0826-3}
\lim _{ t\to 0} 
\Big\|  \psi (2^{-N} \sqrt{A})
   \Big(e^{-t A^{\alpha /2} } f - f  \Big) 
   \Big\|_{L^p}  = 0
   \quad \text{for each } N \in \mathbb N.
\end{equation}
Indeed, we take $\ell_0 \in \mathbb N$ with $2\ell_0  > \alpha$, and write by the subordination 
\[
\begin{split}
\psi (2^{-N} \sqrt{A})
   \Big(e^{-t A^{\alpha /2} } f - f  \Big)  
   =
&    \int_0 ^\infty F_{t,\frac{\alpha}{2 \ell_0}} (s) 
     \psi (2^{-N} \sqrt{A}) \big( e^{-s A^{ \ell_0}} -1 \big)f \,{\rm d}s
\\[1ex]
   =
&    \int_0 ^\infty F_{1,\frac{\alpha}{2 \ell_0}} (s) 
     \psi (2^{-N} \sqrt{A}) \big( e^{-s t^{1/\alpha} A^{ \ell_0}} -1 \big)f \,{\rm d}s .
\end{split}
\]
%
%
%
By extending $\psi (2^{-N}\sqrt{\lambda})$ $(\lambda \geq 0)$ 
to a function in $C_0^\infty (\mathbb R)$,  we have that 
\[
 \| (1+\lambda^2)^\beta \psi(2^{-N}\sqrt{\lambda})\big( e^{-st^{1/\alpha} \lambda^{\ell_0}}-1\big)\|_ {H^\gamma(\lambda\in\mathbb R)} \to 0
\]
as $t \to 0$, where $\beta , \gamma > 1$. Then the 
boundedness of the spectral multiplier in Lemma~\ref{lem:1} implies 
that 
\[
\Big\|  \psi (2^{-N} \sqrt{A}) \big( e^{-s t^{1/\alpha} A^{ \ell_0}} -1 \big) 
\Big\|_{L^p \to L^p} 
   \to 0
\]
as $t\to 0$, and we obtain \eqref{0826-3}. Therefore, the case when $1 \leq p < \infty$ is proved. 

The case when $p = \infty$ follows by duality.  In fact, let  
$f \in L^\infty (\Omega)$ and $g \in L^1 (\Omega)$. Then 
\begin{align*}
\Big| \int_{\Omega} 
\big( e^{-t A^{\frac{\alpha}{2}}} f - f \big) \overline{ g} \,{\rm d}x
\Big| 
& = \Big| \int_{\Omega} 
f \cdot  \overline{ \big( e^{-t A^{\frac{\alpha}{2}}} g - g \big) } \,{\rm d}x \Big|
\\
& \leq  \| f \|_{L^\infty} 
  \| e^{-t A^{\frac{\alpha}{2}}}  g - g \|_{L^1} \\
& \to 0 \quad \text{as } t \to 0 .
\end{align*}
Concerning the above duality argument note that for a sequence $\{g_n\} \subset C_0^\infty (\Omega)$ satisfying $g_n \to g$ in $L^1 (\Omega)$ there holds
\[
\int_{\Omega} 
\big( e^{-t A^{\frac{\alpha}{2}}} f \big)  \overline{ g_n} \,{\rm d}x
= \,_{\mathcal X'}\big\langle e^{-t A^{\frac{\alpha}{2}}} f , g_n \big\rangle_{\mathcal X}
= \,_{\mathcal X'}\big\langle  f , e^{-t A^{\frac{\alpha}{2}}}g_n \big\rangle_{\mathcal X}
= \int_{\Omega} 
 f  \cdot \overline{ e^{-t A^{\frac{\alpha}{2}}} g_n} \,{\rm d}x.
\]
As $n\to\infty$, we obtain that
$\int_{\Omega} 
\big( e^{-t A^{\frac{\alpha}{2}}} f  \big)  \overline{ g} \,{\rm d}x
=  \int_{\Omega} 
 f  \cdot \overline{ e^{-t A^{\frac{\alpha}{2}}} g} \,{\rm d}x $.
\hfill $\Box$ 

\vskip3mm

\noindent 
{\bf Proof of Theorem~\ref{thm:3} (iv), homogeneous case. }  
We start by proving that 
\begin{equation}\label{0927-7}
\sup _{\lambda \in \mathbb C , \, {\rm Re} \,\lambda > 0}\| \lambda (\lambda + A)^{-1} \|_{L^p \to L^p} < \infty. 
\end{equation}
Let $f \in L^p (\Omega)$. 
By the resolution of the identity (see Lemma~\ref{lem:5}) 
for the operator $A/{\rm Re} \lambda$, 
\[
f = \psi ( \sqrt{A/{\rm Re} \lambda}) f 
+ \sum _{j=1}^\infty \phi_j (\sqrt{A/{\rm Re} \lambda}) f .
\]
We note that this holds not only in $\mathcal X'(A)$ but also pointwise almost everywhere, 
and 
\[
\lambda (\lambda + A)^{-1} 
= \frac{\lambda / {\rm Re } \lambda}{ \lambda / {\rm Re} \lambda + A/{\rm Re } \lambda}  
\]
with an abuse of notation. By the spectral multiplier theorem, 
if ${\rm Re} \, \lambda \geq|{\rm Im } \, \lambda|$, \RED  \BLACK 
\[
\Big\| \frac{\lambda / {\rm Re } \lambda}{ \lambda / {\rm Re} \lambda + A/{\rm Re } \lambda}
   \psi ( \sqrt{A/{\rm Re} \lambda}) f 
\Big\|_{L^p} 
\leq C \| f \|_{L^p}, 
\]
\[
\sum_{j=1}^\infty 
\Big\|  \frac{\lambda / {\rm Re } \lambda}{ \lambda / {\rm Re} \lambda + A/{\rm Re } \lambda}
\phi_j (\sqrt{A/{\rm Re} \lambda}) f
\Big\| _{L^p} 
\leq \sum_{j = 1}^\infty \frac{C}{2^{2j}} \| f \|_{L^p}
\leq C \| f \|_{L^p} ,
\]
which proves \eqref{0927-7} on a sector in the right half plane of opening angle $\pi/2$ and hence the sectoriality of $A$ from $ \dot H^s_p (A) \cap \dot H^{s+2}_p (A) $ to $ \dot H^s_p (A)$. 
If ${\rm Re} \, \lambda \leq |{\rm Im } \, \lambda|$,  then we can obtain a similar inequality, 
considering $|{\rm Im } \, \lambda|$ instead of ${\rm Re} \, \lambda$. 
This proves the analyticity of $\{ e^{-tA} \}_{t \geq 0}$ on $(0,\infty)$, see \cite[Chapter 2]{Lunardi_1995}. 
Due to the spectral multiplier estimates in Lemma~\ref{lem:1}, this argument applies to $A^{\alpha /2}$ 
from $ \dot H^s_p (A) \cap \dot H^{s+\alpha}_p (A) $ to $ \dot H^s_p (A)$, 
as well.  
\hfill $\Box$

\vskip3mm

\noindent 
{\bf Proof of Theorem~\ref{thm:3} (i) - (iv), non-homogeneous case. } 
The proof follows
the homogeneous case above except for the local-in-time estimate (ii). 
This is because of the low frequency estimate; in fact, 
we estimate by Lemma~\ref{lem:0912-2} and Lemma \ref{lem:1}
\[
\begin{split}
\| e^{-tA^{\alpha /2}}f \| _{H^{s_2}_{p_2}(A)}
\leq 
& 
C\| \psi (\sqrt{A}) e^{-tA^{\alpha /2}}f \| _{L^{p_2}}
+ C\Big\| \sum _{ j=1}^\infty \phi_j (\sqrt{A}) e^{-tA^{\alpha /2}} f \Big\|_{\dot H^{s_2}_{p_2}(A)}
\\
\leq 
& C t^{-\frac{d}{\alpha} (\frac{1}{p_1} - \frac{1}{p_2})} 
   \| \psi (\sqrt{A}) f \|_{L^{p_1}} 
  + C t^{ - \frac{s_2 -s_1}{\alpha}-\frac{d}{\alpha} (\frac{1}{p_1} - \frac{1}{p_2})} 
    \Big\| \sum _{j=1}^\infty \phi_j (\sqrt{A}) f \Big\| _{\dot H^{s_1}_{p_1}(A)}
\\
\leq 
& C 
\Big( 1 + t^{- \frac{s_2 - s_1}{\alpha}} \Big) 
t^{-\frac{d}{\alpha} (\frac{1}{p_1} - \frac{1}{p_2}) 
         }      \| f \|_{ H^{s_1}_{p_1} (A) }  .
\end{split}
\]
\hfill $\Box$

\vskip3mm 

\noindent
{\bf Acknowledgements. }
T. Iwabuchi was supported by the Grant-in-Aid for Young Scientists (A) (No.~17H04824)
from JSPS.


\begin{bibdiv}
\begin{biblist}

\bib{AS-1961}{article}{
   author={Aronszajn, N.},
   author={Smith, K. T.},
   title={Theory of Bessel potentials. I},
   journal={Ann. Inst. Fourier (Grenoble)},
   volume={11},
   date={1961},
   pages={385--475},
}

\bib{BrMi-2018}{article}{
   author={Brezis, Ha\"{\i}m},
   author={Mironescu, Petru},
   title={Gagliardo-Nirenberg inequalities and non-inequalities: the full
   story},
   journal={Ann. Inst. H. Poincar\'{e} Anal. Non Lin\'{e}aire},
   volume={35},
   date={2018},
   number={5},
   pages={1355--1376},
}

\bib{BrMi-2019}{article}{
   author={Brezis, Ha\"{\i}m},
   author={Mironescu, Petru},
   title={Where Sobolev interacts with Gagliardo--Nirenberg},
   journal={J. Funct. Anal.},
   volume={277},
   date={2019},
   number={8},
   pages={2839--2864},
}

\bib{Ca-1961}{article}{
   author={Calder\'{o}n, A.-P.},
   title={Lebesgue spaces of differentiable functions and distributions},
   conference={
      title={Proc. Sympos. Pure Math., Vol. IV},
   },
   book={
      publisher={Amer. Math. Soc., Providence, R.I.},
   },
   date={1961},
   pages={33--49},
}

\bib{Chi-2018}{article}{
   author={Chikami, Noboru},
   title={On Gagliardo-Nirenberg type inequalities in Fourier-Herz spaces},
   journal={J. Funct. Anal.},
   volume={275},
   date={2018},
   number={5},
   pages={1138--1172},
}

\bib{Davies_1990}{book}{
   author={Davies, E. B.},
   title={Heat kernels and spectral theory},
   series={Cambridge Tracts in Mathematics},
   volume={92},
   publisher={Cambridge University Press, Cambridge},
   date={1990},
}

\bib{DiPaVa-2012}{article}{
   author={Di Nezza, Eleonora},
   author={Palatucci, Giampiero},
   author={Valdinoci, Enrico},
   title={Hitchhiker's guide to the fractional Sobolev spaces},
   journal={Bull. Sci. Math.},
   volume={136},
   date={2012},
   number={5},
   pages={521--573},
}

\bib{Gag-1959}{article}{
   author={Gagliardo, Emilio},
   title={Ulteriori propriet\`a di alcune classi di funzioni in pi\`u variabili},
   language={Italian},
   journal={Ricerche Mat.},
   volume={8},
   date={1959},
   pages={24--51},
}

\bib{HMOW-2011}{article}{
   author={Hajaiej, Hichem},
   author={Molinet, Luc},
   author={Ozawa, Tohru},
   author={Wang, Baoxiang},
   title={Necessary and sufficient conditions for the fractional
   Gagliardo-Nirenberg inequalities and applications to Navier-Stokes and
   generalized boson equations},
   conference={
      title={Harmonic analysis and nonlinear partial differential equations},
   },
   book={
      series={RIMS K\^{o}ky\^{u}roku Bessatsu, B26},
      publisher={Res. Inst. Math. Sci. (RIMS), Kyoto},
   },
   date={2011},
   pages={159--175},
}

\bib{Iwa-2018}{article}{
   author={Iwabuchi, Tsukasa},
   title={The semigroup generated by the Dirichlet Laplacian of fractional
   order},
   journal={Anal. PDE},
   volume={11},
   date={2018},
   number={3},
   pages={683--703},
}


\bib{IMT-2018}{article}{
   author={Iwabuchi, Tsukasa},
   author={Matsuyama, Tokio},
   author={Taniguchi, Koichi},
   title={Boundedness of spectral multipliers for Schr\"{o}dinger operators on
   open sets},
   journal={Rev. Mat. Iberoam.},
   volume={34},
   date={2018},
   number={3},
   pages={1277--1322},
}

\bib{IMT-pre}{article}{
   author={Iwabuchi, Tsukasa},
   author={Matsuyama, Tokio},
   author={Taniguchi, Koichi},
   title={$L^p$-mapping properties for the Schr\"odinger operators in open sets of $\mathbb R^d$},
   journal={preprint, arXiv:1602.08208},
}

\bib{IMT-2019}{article}{
   author={Iwabuchi, Tsukasa},
   author={Matsuyama, Tokio},
   author={Taniguchi, Koichi},
   title={Besov spaces on open sets},
   journal={Bull. Sci. Math.},
   volume={152},
   date={2019},
   pages={93--149},
}

\bib{Lunardi_1995}{book}{
   author={Lunardi, Alessandra},
   title={Analytic semigroups and optimal regularity in parabolic problems},
   series={Modern Birkh\"{a}user Classics},
   note={[2013 reprint] },
   publisher={Birkh\"{a}user/Springer Basel AG, Basel},
   date={1995},
}

\bib{MaOz-2003}{article}{
   author={Machihara, Shuji},
   author={Ozawa, Tohru},
   title={Interpolation inequalities in Besov spaces},
   journal={Proc. Amer. Math. Soc.},
   volume={131},
   date={2003},
   number={5},
   pages={1553--1556},
}

\bib{Nir-1959}{article}{
   author={Nirenberg, L.},
   title={On elliptic partial differential equations},
   journal={Ann. Scuola Norm. Sup. Pisa Cl. Sci. (3)},
   volume={13},
   date={1959},
   pages={115--162},
}


\bib{Sobo-1935}{article}{
   author={Sobolev, Sergei Lvovich},
	 title={Le syst\`eme de Cauchy dans l'espace des fonctionelles},
   journal={Dokl. Akad. Nauk SSSR},
   volume={3},
   date={1935},
   pages={291--294},
}

\bib{Sobo-1936}{article}{
   author={Sobolev, Sergei Lvovich},
   title={M\'ethode nouvelle \`a r\'esoudre le probl\`eme de Cauchy pour les \'equations lin\'eaires 
   hyperboliques normales},
   journal={Mat. Sb.},
   volume={1},
   date={1936},
   pages={39--72},
}

\bib{Sobo-1938}{article}{
   author={Sobolev, Sergei Lvovich},
   title={Sur un th\'eor\`eme d'analyse fonctionelle },
   journal={Mat. Sb.},
   volume={4},
   date={1938},
   pages={471--497},
}


\bib{Ryc-1999}{article}{
    AUTHOR = {Rychkov, V. S.},
     TITLE = {On restrictions and extensions of the {B}esov and
              {T}riebel-{L}izorkin spaces with respect to {L}ipschitz
              domains},
   JOURNAL = {J. London Math. Soc. (2)},
    VOLUME = {60},
      YEAR = {1999},
     PAGES = {237--257},
} 

\bib{TOS-2002}{article}{
   author={Thinh Duong, Xuan},
   author={Ouhabaz, El Maati},
   author={Sikora, Adam},
   title={Plancherel-type estimates and sharp spectral multipliers},
   journal={J. Funct. Anal.},
   volume={196},
   date={2002},
   number={2},
   pages={443--485},
}

\bib{Triebel_1983}{book}{
   author={Triebel, Hans},
   title={Theory of function spaces},
   series={Monographs in Mathematics},
   volume={78},
   publisher={Birkh\"auser Verlag, Basel},
   date={1983},
}

\bib{Triebel_1995}{book}{
   author={Triebel, Hans},
   title={Interpolation theory, function spaces, differential operators},
   edition={2},
   publisher={Johann Ambrosius Barth, Heidelberg},
   date={1995},
   pages={532},
}


\bib{Yosida_1980}{book}{
   author={Yosida, K\^{o}saku},
   title={Functional analysis},
   series={Grundlehren der Mathematischen Wissenschaften [Fundamental
   Principles of Mathematical Sciences]},
   volume={123},
   edition={6},
   publisher={Springer-Verlag, Berlin-New York},
   date={1980},
}

\end{biblist}
\end{bibdiv}
\end{document}